\documentclass[english]{article}

\usepackage[T1]{fontenc}
\usepackage[latin9]{inputenc}
\usepackage{geometry}
\usepackage{color}
\usepackage{amsmath}
\usepackage{amssymb}
\usepackage{esint}
\usepackage{helvet}

\makeatletter
\usepackage{amsthm}\usepackage{amsfonts}

\newtheorem{theorem}{Theorem}[section]

\newtheorem{corollary}{Corollary}[section]
\newtheorem{remark}{Remark}[section]
\usepackage{babel}

\theoremstyle{definition}
\newtheorem{example}{Example}[section]

\makeatother

\begin{document}

\thispagestyle{empty}

\title{Estimates for  Kantorovich functionals between solutions to Fokker
-- Planck -- Kolmogorov equations with dissipative drifts}

\author{Oxana A. Manita}
\maketitle
\begin{abstract}
We obtain estimates for the Kantorovich functionals between solutions
to different Fokker -- Planck -- Kolmogorov equations for measures
with same diffusion part but different drifts and different initial
conditions. We show possible applications of such estimates to the
study of the well-posedness for nonlinear equations. 
\end{abstract}

\textsc{Keywords:} Fokker -- Planck -- Kolmogorov equation; Kantorovich distance;
 Nonlinear Fokker -- Planck equation;  Dissipative operator.

\textsc{Author's address:} Lomonosov Moscow State University, Faculty
of Mechanics and Mathematics, Russia, 119991, Moscow, GSP-1, 1, Leninskiye
Gory, Main Building; o.manita@lambda.msu.ru.

\section{ Introduction.}

In the present paper we derive and study estimates for the Kantorovich
functionals between probability solutions for the linear 
Fokker -- Planck -- Kolmogorov (FPK) equations for probability measures $\mu_t$ and $\sigma_t$ on $\mathbb{R}^d$,
 $t\in[0,T]$, 
 with 
different drifts and different initial conditions 
\begin{eqnarray}
\partial_{t}\mu_{t} & = & \mbox{trace}(Q(x,t)D^{2}\mu_{t})-\mbox{div}(B_{1}(x,t)\mu_{t}),
\qquad\mu|_{t=0}=\mu_{0}\nonumber\\
\partial_{t}\sigma_{t} & = & \mbox{trace}(Q(x,t)D^{2}\sigma_{t})-\mbox{div}(B_{2}(x,t)\sigma_{t}),
\qquad\sigma|_{t=0}=\sigma_{0}.
\nonumber\end{eqnarray}
 We also show an
alternative method to the study of well-posedness and stability of
solutions to the nonlinear FPK equations
\begin{equation}
\partial_{t}\rho_{t}=\mbox{trace}(Q(x,t)D^{2}\rho_{t})-\mbox{div}(B(\rho,x,t)\rho_{t}),
\quad\rho|_{t=0}=\rho_{0},
\label{eq:e1}\end{equation}
based on such estimates
for linear equations.

Recently, FPK equations have been actively studied from the functional-analytical,
variational and as well from the probabilistic point of view. Interesting
connections between approaches have been found (a survey of the current
state of studies is provided in \cite{FPbook}). Estimates connecting
distances between solutions with distances between initial data and
even coefficients play a great role not only for the study of such
qualitative properties of solutions as uniqueness or stability, but
also for numerical simulations. In this context estimates for distances
between solutions to equations with different drift terms are particularly
interesting.

In Section 1 we derive estimates for the Kantorovich functionals between
solutions of FPK equations with different dissipative drifts. To do
this, we partially use ideas from \cite{NPS}. Since these ideas can
not be directly applied neither in the case of different drifts nor
to nonlinear equations, new methods and ideas should be used. Extension
to these cases has been done for the Kantorovich functionals with
bounded cost functions. Moreover, we admit time-dependent coefficients
and a non-unit diffusion matrix $Q$. We note that the requirement
of dissipativity is not really restrictive -- in most physical examples,
the drift term is a minus  gradient of a convex function, i.e. dissipative.
Section 2 is concerned with applications of these estimates to the
study of the well-posedness of the Cauchy problem for the nonlinear
FPK equation. Well-posedness for the nonlinear equations has been
studied by many authors even in a more general setting (see, for example,
\cite{D,ManShap,Funaki,ManRomShap}). However we present an alternative
approach to this problem that is applicable in case of dissipative
drifts. A similar method of treating well-posedness via estimates
for the distances between solutions to linear equations was used in
\cite{BSh_estimates}.

Let us introduce some notation and give basic definitions. By $C_{0}^{\infty}(\mathbb{R}^{d})$
and $C_{0}^{\infty}(\mathbb{R}^{d}\times(0,T))$ we denote
classes of infinitely smooth compactly supported functions on $\mathbb{R}^{d}$
and $\mathbb{R}^{d}\times(0,T)$ respectively. For shortness
of notation we shall always drop
the subscript $\mathbb{R}^{d}$ when integrating over the whole space. 
We shall say that a measure $\rho$ on $\mathbb{R}^{d}\times[0,T]$
is given by a family of probability measures $(\rho_{t})_{t\in[0,T]}$
on $\mathbb{R}^{d}$ (and write $\rho(dx\, dt)=\rho_{t}(dx)\, dt$
or simply $\rho=\rho_{t}dt$), if $\rho_{t}\ge0$, $\rho_{t}(\mathbb{R}^{d})=1$,
for each Borel set $U$ the function $t\mapsto\rho_{t}(U)$ is measurable
and 
$$
\int_{0}^{T}\int\phi\, d\rho=\int_{0}^{T}\int\phi\, d\rho_{t}\, dt\quad
\forall\phi\in C_{0}^{\infty}(\mathbb{R}^{d}\times(0,T)).
$$
Given a probability measure $\rho_{0}$ on $\mathbb{R}^{d}$, a symmetric Borel
matrix $Q(x,t)$ and a Borel mapping $B(x,t):\,\mathbb{R}^{d}\times[0,T]\rightarrow\mathbb{R}^{d}$,
consider the following Cauchy problem for the linear FPK equation
\begin{equation}
\partial_{t}\rho_{t}=\mbox{trace}(Q(x,t)D^{2}\rho_{t})-\mbox{div}(B(x,t)\rho_{t}),\quad\rho|_{t=0}=\rho_{0}.
\label{eq:lin1}\end{equation}
Here $D^{2}$ denotes the Hessian matrix with respect to the spacial
variables. Denote the elements of the diffusion matrix $Q(x,t)$ by
$q^{ij}(x,t),\,1\leq i,j\leq d$ and the elements of the vector drift
$B(x,t)$ by $b^{j}(x,t),\,1\leq j\leq d$. Set 
$$
L\phi=q^{ij}(x,t)\partial_{x_{i}x_{j}}^{2}\phi+b^{i}(x,t)\partial_{x_{i}}\phi,
$$
where summation over all repeated indices is taken. We shall say that
a measure $\rho(dx\, dt)=\rho_{t}(dx)\, dt$ is a solution to the
Cauchy problem (\ref{eq:lin1}), if the mappings $q^{ij}(x,t),\,\, b^{i}(x,t),\,1\leq i,j\leq d,$
are Borel and belong to $L^{1}(\rho,U\times[0,T])$ for each ball
$U\subset\mathbb{R}^{d}$, and for each test function $\varphi\in C_{0}^{\infty}(\mathbb{R}^{d})$
we have 
\begin{equation}
\int\varphi\, d\rho_{t}=\int\varphi\, d\rho_{0}+\int_{0}^{t}\int L\varphi\, d\rho_{s}\, ds
\label{r1}\end{equation}
 for all $t\in[0,T]$. Sometimes it is more convenient to use an equivalent
definition (see \cite{obzor-Uni}), more precisely, the identity 
\begin{equation}
\int\phi(x,t)\, d\rho_{t}=\int\phi(x,0)\, d\rho_{0}+
\int_{0}^{t}\int\bigl[\partial_{s}\phi+L\phi\bigr]\, d\rho_{s}\, ds,
\label{r2}\end{equation}
 for all $t\in[0,T]$ and all test functions 
 $\psi\in C^{2,1}(\mathbb{R}^{d}\times[0,T))\bigcap C(\mathbb{R}^{d}\times[0,T])$
that are identically zero outside some ball $U\subset\mathbb{R}^{d}$.
If we know a priori that the drift term $B$ is integrable over $\mathbb{R}^{d}\times[0,T]$
with respect to the measure $d\rho_{s}ds$ and $\phi$ is supported
on the whole $\mathbb{R}^{d}$, but has two continuous bounded derivatives, then
\eqref{r2} also holds true for such $\phi$ (to show this, it sufficies
to use a standard truncation argument).

\section{Estimates for the Kantorovich functionals between solutions to linear
equations with different drifts}

In this section, we shall focus on two solutions of the linear FPK
equation with different initial conditions and different drifts. Fix
$T>0$. Given probability measures $\mu_{0}$ and $\sigma_{0}$ on
$\mathbb{R}^{d}$, a symmetric Borel matrix $Q(x,t)$ and Borel mappings $B_{\mu},\, B_{\sigma}:\mathbb{R}^{d}\times[0,T]\rightarrow\mathbb{R}^{d}$,
consider two Cauchy problems 
\begin{eqnarray}
\partial_{t}\mu_{t} & = & \mbox{trace}(Q(x,t)D^{2}\mu_{t})-\mbox{div}(B_{\mu}(x,t)\mu_{t}),
\qquad\mu|_{t=0}=\mu_{0}\label{eq:2lin}\\
\partial_{t}\sigma_{t} & = & \mbox{trace}(Q(x,t)D^{2}\sigma_{t})-\mbox{div}(B_{\sigma}(x,t)\sigma_{t}),
\qquad\sigma|_{t=0}=\sigma_{0}.
\nonumber \end{eqnarray}
We emphasize that the indices $\mu$ and $\sigma$ in the drift coefficients 
are merely used to distinguish the different drifts (by marking corresponding
solutions), and it is not necessary to define $B$ as a map on a space
of measures.

Given a monotone nonnegative continuous function $h$ on $\mathbb{R}$
with $h(0)=0$, introduce the Kantorovich $h$-cost functional between
the probability measures $\mu$ and $\sigma$ by 
\begin{equation}
C_h (\mu,\sigma):=\inf_{\pi\in\Pi(\mu,\sigma)}\int_{\mathbb{R}^{d}\times\mathbb{R}^{d}}h(|x-y|)d\pi(x,y),
\label{eq:ch-1}\end{equation}
 where $\Pi(\mu,\sigma)$ is the set of couplings between
$\mu$ and $\sigma$. Recall that a probability measures $\pi$ on
$\mathbb{R}^{d}\times\mathbb{R}^{d}$ belongs to $\Pi(\mu,\sigma)$ iff $\pi(E\times\mathbb{R}^{d})=\mu(E)$,
$\pi(\mathbb{R}^{d}\times E)=\sigma(E)$ for each Borel
set $E\subset\mathbb{R}^{d}$. If $h$ is a concave function
with $h(r)>0$ for $r>0$, then $C_{h}$ defines a distance on the
space of probability measures and turns it into a complete metric space with
topology that coincides with the usual weak one (see \cite[Proposition 7.1.5]{AGS}).
Another important example is given by $h(r)=\min\{ |r|^{p},1\} $
for some $p\geq1$. In this case $C_{h}^{1/p}$ turns the space of
probability measures into a complete metric space. Moreover, convergence
with respect to this metric is equivalent to the weak convergence (see
\cite[Th. 1.1.9]{obzor_Monge}). 

Further we assume that a monotone non-decreasing continuous bounded
cost function $h$ with $h(0)=0$ is fixed. Set $\Vert h\Vert _{\infty}:=\sup_{z\in\mathbb{R}^{d}}h(|z|)<\infty$.

Throughout the paper we assume that the following regularity condition
holds:

$({\rm A1})$ The diffusion matrix $Q(x,t)$ has uniformly
bounded elements with uniformly bounded first derivatives. Moreover,
it is strictly elliptic: there exists $\nu>0$ such that  $\forall(x,t)\in\mathbb{R}^{d}\times[0,T]$ 
\begin{equation}
\langle Q(x,t)y,y\rangle \geq\nu|y|^{2}\quad\forall y\in\mathbb{R}^{d}.
\label{eq:elliptic}\end{equation}

\begin{theorem}\label{est} Let ${\rm (A1)}$ hold. Let $(\mu_{t})_{t\in[0,T]}$
and $(\sigma_{t})_{t\in[0,T]}$ be solutions
to \eqref{eq:2lin} with initial conditions $\mu_{0}$ and $\sigma_{0}$
respectively. Suppose that the drift term $B_{\mu}$ is $\lambda$-dissipative
in $x$, i.e. 
\begin{equation}
\langle B_{\mu}(x,t)-B_{\mu}(y,t),x-y\rangle \leq\lambda\Vert x-y\Vert ^{2}
\label{eq:B}\end{equation}
 for all $x,y\in\mathbb{R}^{d}$ and all $t\in[0,T]$. Let 
 \begin{equation}
B_{\mu}(x,t)-\lambda x,\,\, B_{\sigma}(x,t)-\lambda x\in L^{2}(\mathbb{R}^{d}\times[0,T],d(\mu_{s}+\sigma_{s})ds)
\label{eq:B-reg-1}\end{equation}
 Then 
 \begin{equation}
C_{h_{\lambda t}} (\mu_t,\sigma_t) \leq C_h (\mu_0,\sigma_0)
+\Vert h\Vert _{\infty}\cdot\sqrt{\int_{0}^{t}\int\nu^{-1}|B_{\mu}-B_{\sigma}|^{2}d\sigma_{s}ds}
\sqrt{1+\int_{0}^{t}\int\nu^{-1}|B_{\mu}-B_{\sigma}|^{2}d\sigma_{s}ds},
\label{eq:basic_est}\end{equation}
for all $t\in[0,T]$, where ${\displaystyle h_{s}(r):=h(re^{-s})}$.

\end{theorem}

\begin{remark} The bound \eqref{eq:basic_est} is obviously asymmetric
in measure:  we impose dissipativity on $B_{\mu}$, and the integration
in the right-hand side is taken over $\sigma$. This property might
be interesting from the point of view of possible numerical simulations.
Indeed, if we want to solve a FPK equation 
$$
\partial_{t}\mu_{t}={\rm \mbox{trace}}(Q(x,t)D^{2}\mu_{t})-{\rm \mbox{div}}(B(x,t)\mu_{t})
$$
with a dissipative drift $B$, we can approximate the drift with {}``better''
drifts $B_{n}$ and solve FPK equations with those drifts. Then \eqref{eq:basic_est}
controls the distance between the desired solution $\mu_{t}$ and
the approximative solution $\mu_{t}^{n}$ in terms of the distance
between drifts integrated over the \textbf{known} solution
$\mu^{n}$. \end{remark}

\textbf{Proof.} Let $(\mu_{t})_{t\in[0,T]}$ and $(\sigma_{t})_{t\in[0,T]}$
satisfy the assumptions of the theorem. By virtue of \cite{BKR} measures
$\mu_{t}$ and $\sigma_{t}$ have strictly positive densities with
respect to Lebesgue measure on $\mathbb{R}^{d}$ for each $t\in[0,T]$. We split
the proof of \eqref{eq:basic_est} into several steps.

\textbf{Step 1. Reduction to the dissipative case ($\lambda=0$).}
We rescale the problem, keeping the cost function unchanged, in order
to reduce the problem to the case of a dissipative drift $B_{\mu}$.
To this aim we use the rescaling procedure from \cite{NPS} with the
opposite sign (since our drift term and the drift term in the cited
work have the opposive signs). For completeness, we provide the rescaling
procedure: for $\lambda\neq0$ define the change of time 
$$
s(t):=\int_{0}^{t}e^{-2\lambda r}dr=\frac{1-e^{-2\lambda t}}{2\lambda},\quad t(s)=\frac{-\ln(1-2\lambda s)}{2\lambda},\,\, s\in[0,S_{\infty}),
$$
where $S_{\infty}=+\infty$ for $\lambda<0$ and $S_{\infty}=1/(2\lambda)$
for $\lambda>0$. For measures $\mu_{t}$ and $\sigma_{t}$
introduce their rescaled versions $\rho_{s}^{\mu}$ and $\rho_{s}^{\sigma}$:
for each Borel set $E\subset\mathbb{R}^{d}$ define 
${\displaystyle \rho_{s}^{w}(E):=w_{t(s)}(e^{\lambda t(s)}E)}$
for $w=\mu,\sigma$. 
We notice that 
${\displaystyle C_{h} (\rho^{\mu}_s,\rho^{\sigma}_s)= C_{h_{\lambda t}} (\mu_t,\sigma_t).}$
Since $B_{\mu}$ is $\lambda$-dissipative, $A_{\mu}:=B_{\mu}-\lambda I$
is dissipative. Define the rescaled diffusion coefficient by 
$$
\tilde{Q}(y,s):=Q(t(s),e^{\lambda t(s)}y)
$$
and the rescaled drifts by 
$$
\tilde{B}_{w}(y,s):=e^{\lambda t(s)}B_{w}(t(s),e^{\lambda t(s)}y),\,\,\tilde{A}_{w}(y,s):=e^{\lambda t(s)}A_{w}(t(s),e^{\lambda t(s)}y),\, w=\mu,\sigma.
$$
Note that $\tilde{A}_{\mu}$ is also a dissipative operator. The
measure $\mu=\mu_{t}dt$ is a solution to 
$$
\partial_{t}\mu_{t}=\mbox{trace}(QD^{2}\mu_{t})-\mbox{div}(B_{\mu}\mu_{t})
$$
iff the rescaled version $\rho^{\mu}=\rho_{t}^{\mu}dt$ is a solution
to 
\begin{equation}
\partial_{t}\rho_{t}^{\mu}=\mbox{trace}(\tilde{Q}D^{2}\rho_{t}^{\mu})-\mbox{div}(\tilde{A}_{\mu}\rho_{t}^{\mu});
\label{eq:tozh_resc}\end{equation}
moreover, \eqref{eq:B-reg-1} holds true iff for all nonnegative
$s_{1}<s_{2}\leq S(T)<S_{\infty}$ one has 
$$
\int_{s_{1}}^{s_{2}}\int|\tilde{A}_{\mu}(x,s)|^{2}d\rho_{s}^{\mu}ds<+\infty.
$$
The integrability statement follows immediately from the change of
variables formula, the identity \eqref{eq:tozh_resc} can be checked
explicitly: it sufficies to consider the change of variables $\mathtt{X}(x,t):=(e^{-\lambda t}x,s(t))$
and calculate the derivatives. Similar statement holds for $\sigma$
and $\rho^{\sigma}$. This means that it is sufficient to prove \eqref{eq:basic_est}
only in the case $\lambda=0$, i.e. in the case of a dissipative drift
term $B_{\mu}$.

\textbf{Step 2. Approximation of the drift term.} We construct a family
of smooth (in both variables) bounded Lipschitz (as functions of $x$) dissipative
operators $A_{k}^{\varepsilon}(x,t)$, approximating the dissipative
drift term $B_{\mu}(x,t)$. 

For each $t$ the operator $B_{\mu}(\cdot,t)$ can be approximated
by Lipschitz in $x$ bounded dissipative operators $A_{k}(\cdot,t)$
with bounded first order derivatives with respect to the spacial variables (see
\cite[Th. 2.4, 2.5]{NPS}): for each $t\in [0,T]$
\begin{equation}
\lim_{k\rightarrow\infty}A_{k}(x,t)=B_{\mu}(x,t)\,\,\mbox{for a.e.}\,\,(x)\in\mathbb{R}^{d}\quad\mbox{and }\sup_{x\in\mathbb{R}^{d}}|A_{k}(x,t)|\leq k+1.
\label{eq:smooth1}\end{equation}
Fix some non-negative function $\eta\in C_{0}^{\infty}([0,T])$
such that $\Vert \eta\Vert _{L^{1}([0,T])}=1$.
Introduce a family of mollifiers $\eta_{\varepsilon}(t):=\eta(t/\varepsilon)$.
Since for each $k$ the mapping $A_{k}(x,t)$ is bounded,
the mappings $A_{k}^{\varepsilon}(x,t):=\eta_{\varepsilon}(t)*A_{k}(x,t)$ have
bounded derivatives of all orders with respect to $t$ and converge
to $A_{k}(x,t)$ as $\varepsilon\rightarrow0$ for a.e. $(x,t)\in\mathbb{R}^{d}\times[0,T]$.
Notice that $A_{k}^{\varepsilon}$ also have bounded first order derivatives
with respect to  the spacial variables. Moreover, $A_{k}^{\varepsilon}$ are dissipative
in $x$: 
$$
\langle A_{k}^{\varepsilon}(x,t)-A_{k}^{\varepsilon}(y,t),x-y\rangle =\eta_{\varepsilon}(t)*\langle A_{k}(x,t)-A_{k}(y,t),x-y\rangle \leq0.
$$
 Finally, we define operators $\mathcal{L}_{k}^{\varepsilon}$ as follows: for $(x,s)\in\mathbb{R}^{d}\times[0,T]$
$$
\mathcal{L}_{k}^{\varepsilon}[\phi](x,s):=\mbox{trace}(Q(x,s)D^{2}\phi(x,s))+\langle A_{k}^{\varepsilon}(x,s),\nabla_{x}\phi(x,s)\rangle ,\quad\phi(\cdot,s)\in C^{2}(\mathbb{R}^{d}).
$$

\textbf{Step 3. Reduction of the class of test functions.} It is
well-known (for example,  \cite[Th. 1.3]{Villani}) that the problem
\eqref{eq:ch-1} admits a dual formulation: define the class $\Phi_{h}$
as 
$$
\Phi_{h}:=\{ (\phi,\psi)\in L^{1}(\mu)\times L^{1}(\nu):\phi(x)+\psi(y)\le h(|x-y|)\} .
$$
 Hence 
 \begin{equation}
C_h (\mu,\sigma)=\sup_{(\phi,\psi)\in\Phi_{h}}\int\phi\,\mathrm{d}\mu+\int\psi\,\mathrm{d}\sigma.
\label{eq:duality-1}\end{equation}
 An important observation (\cite[Lemma 2.3]{MZK}) is that in the case
of a bounded cost function $h$ the supremum in the dual problem \eqref{eq:duality-1}
may be taken over a smaller class of functions $\Phi_{h}^{\delta}$
for any $\delta>0$, where 
\begin{equation}
\Phi_{h}^{\delta}:=\Phi_{h}\cap C_{0}^{\infty}(\mathbb{R}^{d})
\cap\{ (\phi,\psi):\,\inf\psi>-\delta\,\,\mbox{and\,\ }\sup\psi\leq\Vert h\Vert _{\infty}\}.
\label{eq:phib}\end{equation}
 The proof is based on the fact that functions $\varphi$ and $\psi$
can be shifted by different constants and truncated in such way that
the new pair $(\varphi_{0},\psi_{0})$ is still admissible
and the value in \eqref{eq:duality-1} doesn't decrease: 
$$
\int\varphi_{0}d\mu+\int\psi_{0}d\sigma\geq\int\varphi d\mu+\int\psi d\sigma\,\,\mbox{and }
\inf_{\mathbb{R}^{d}}\psi_{0}=0,\,\,\sup_{\mathbb{R}^{d}}\psi_{0}\leq\Vert h\Vert _{\infty},\,\,
\sup_{\mathbb{R}^{d}}\varphi_{0}\geq0.
$$
 If we want to deal with smooth compactly supported functions, then
the bounds get a bit worse and lead to the class \eqref{eq:phib}.

In the sequel we shall take the supremum in \eqref{eq:duality-1}
over the class $\Phi_{h}^{b}$ of admissible pairs of $C_{0}^{\infty}(\mathbb{R}^{d})$-functions
such that $\Vert \psi\Vert _{\infty}\leq\Vert h\Vert _{\infty}$.

\textbf{Step 4. Adjoint problem.} Fix an admissible pair $(\phi,\psi)\in\Phi_{h}^{b}$.
The smootheness of operators $A_{k}^{\varepsilon}$ imply \cite[Th. 3.2.1]{StrVar}
that the following adjoint problems have solutions $g,f\in C_{b}^{2,1}(\mathbb{R}^{d}\times[0,t])$:
\begin{equation}
\partial_{s}g+\mathcal{L}_{k}^{\varepsilon} g=0,\quad g(\cdot,t)=\phi(\cdot)\qquad\mbox{and }\qquad\partial_{s}f+\mathcal{L}_{k}^{\varepsilon} f=0,\quad f(\cdot,t)=\psi(\cdot)
\label{eq:backw}\end{equation}
First, due to the maximum principle \cite[Th. 3.1.1]{StrVar} 
\begin{equation}
\sup_{\mathbb{R}^{d}\times[0,t]}|g|\leq\sup_{\mathbb{R}^{d}}|\phi|,\sup_{\mathbb{R}^{d}\times[0,t]}|f|\leq\sup_{\mathbb{R}^{d}}|\psi|.
\label{eq:max_pr}\end{equation}

Let us derive bounds for $|\nabla g|$ and $|\nabla f|$.
The method of doing this is inspired by the Bernstein estimates. Denote
for shortness $A_{k}^{\varepsilon}:=(\alpha^{1},\dots,\alpha^{d})$.
Set $v(x,t):=|\nabla g|^{2}+\kappa g^{2}-t$, where $\kappa$
will be chosen below. Explicit computation gives (summation over all
repeated indices is assumed) 
\begin{equation}
(\partial_{s}-\mathcal{L}_{k}^{\varepsilon})v
\overset{\eqref{eq:backw}}{=}2\partial_{x_{k}} g(\partial_{x_{k}} q^{ij}\partial_{x_{i}x_{j}}^{2}g+
\partial_{x_{k}}\alpha^{i}\partial_{x_{i}}g)-2q^{ij}\partial_{x_{k}x_{i}}^{2}g\partial_{x_{k}x_{j}}^{2}g-
2\kappa q^{ij}\partial_{x_{i}}g\partial_{x_{j}} g-1.
\label{eq:lv}\end{equation}
Due to dissipativity, $DA_{k}^{\varepsilon}$ defines a negative quadratic
form and 
$$
\partial_{x_{k}} g \partial_{x_{k}}\alpha^{i}\partial_{x_{i}}g=
(DA_{k}^{\varepsilon}\nabla g,\nabla g)\leq0.
$$
By virtue of this observation, \eqref{eq:elliptic} and the Cauchy
inequality $2ab\leq ca^{2}+c^{-1}b^{2}$ with $c=2\nu$, \eqref{eq:lv}
is dominated by 
$$
\omega c^{-1}|\nabla g|^{2}+c\sum_{i,j}(\partial_{x_{i}x_{j}}^{2}g)-
2\nu\sum_{i,j}(\partial_{x_{i}x_{j}}^{2}g)-2\nu\kappa|\nabla g|^{2}-1=\Omega|\nabla g|^{2}-2\nu\kappa|\nabla g|^{2}-1,
$$
where $\omega:=2\max\{ |\partial_{x_{k}}q^{ij}|\} $ and
$\Omega:=\omega\cdot(2\nu)^{-1}$ depend \emph{only} on
the diffusion matrix and not on the drift. Choosing $\kappa:=\Omega\cdot(2\nu)^{-1}$,
we get 
$$
(\partial_{s}-\mathcal{L}_{k}^{\varepsilon})v\leq-1.
$$
Therefore the maximum principle \cite[Th. 3.1.1]{StrVar} ensures 
$$
\max_{\mathbb{R}^{d}\times[0,t]}|v|\leq\max_{\mathbb{R}^{d}}|v(x,0)|\equiv\max_{\mathbb{R}^{d}}|\nabla\phi|^{2}+\kappa\max_{\mathbb{R}^{d}}|\phi|^{2},
$$
hence 
\begin{equation}
\sup_{\mathbb{R}^{d}\times[0,t]}|\nabla g(x,s)|\leq(\max_{\mathbb{R}^{d}}|\nabla\phi|^{2}+\kappa\max_{\mathbb{R}^{d}}|\phi|^{2})^{1/2}=:C_{1}.
\label{eq:max_pr_pr_1}\end{equation}
Similarly 
\begin{equation}
\sup_{\mathbb{R}^{d}\times[0,t]}|\nabla f(x,s)|\leq(\max_{\mathbb{R}^{d}}|\nabla\psi|^{2}+\kappa\max_{\mathbb{R}^{d}}|\psi|^{2})^{1/2}+t=:C_{2}.
\label{eq:max_pr_pr}\end{equation}
Set $l:=C_{1}+C_{2}$.

Finally, let us prove the crucial assertion: if the pair $(\varphi,\psi)$
is admissible, and $f$ and $g$ solve \eqref{eq:backw}, then $g(x,0)+f(y,0)\leq h(|x-y|)$.
In the case $Q\equiv I$ it was proved in \cite[Th.3.1]{NPS}. In
the general case the proof almost repeats the case $Q=I$, but we sketch it for completeness.
By approximating $h$ from above, we can assume without lack of generality
that $h\in C^{1}(\mathbb{R})$. Define $H(y_{1},y_{2}):=h(|y_{1}-y_{2}|)$ and 
$$
0\leq\xi(y_{1},y_{2})=\xi(y_{2},y_{1})=\{ \begin{array}{ll}
\frac{h(|y_{1}-y_{2}|)}{|y_{1}-y_{2}|} & \mbox{if }y_{1}\neq y_{2}\\
0 & \mbox{if }y_{1}=y_{2}\end{array}.
$$
First assume that 
$$
\partial_{s}g+\mathcal{L}_{k}^{\varepsilon} g>0,\quad\partial_{s}f+\mathcal{L}_{k}^{\varepsilon} f>0.
$$
Suppose that $\zeta(y_{1},y_{2},s):=g(y_{1},s)+f(y_{2},s)-H(y_{1},y_{2})$
attains a local maximum at $(Y_{1},Y_{2},S)$ and $S<t$. Then 
$\partial_{s}\zeta(Y_{1},Y_{2},S)=\partial_{s}g(Y_{1},S)+\partial_{s}f(Y_{2},S)\leq0$,
$$
\nabla_{y_{1}}\zeta(Y_{1},Y_{2},S)=\nabla_{y_{2}}\zeta(Y_{1},Y_{2},S)=0\Longrightarrow
\nabla_{y_{1}}g(Y_{1},S)=-\nabla_{y_{2}}f(Y_{2},S)=\xi(Y_{1},Y_{2})(Y_{1}-Y_{2})
$$
 and, due to dissipativity, 
$$
A_{k}^{\varepsilon}(Y_{1},S)\nabla_{y_{1}}g(Y_{1},S)+A_{k}^{\varepsilon}(Y_{2},S)\nabla_{y_{2}}f(Y_{2},S)
=\xi(Y_{1},Y_{2})\langle A_{k}^{\varepsilon}(Y_{1},S)-A_{k}^{\varepsilon}(Y_{2},S),Y_{1}-Y_{2}\rangle \leq0.
$$
Since $\zeta(Y_{1}+z,Y_{2}+z,S)$ as a function of $z$ has a local
maximum at $z=0$ and $Q$ is positive definite, 
$\mbox{trace}\bar{Q}D^{2}\zeta=\mbox{trace}Q(Y_{1},S)D^{2}g+\mbox{trace}Q(Y_{2},S)D^{2}f\leq0$,
where 
$$
\bar{Q}(y_{1},y_{2},s):=\Bigl(\begin{array}{ll}
Q(y_{1}) & 0\\
0 & Q(y_{2})\end{array}\Bigr).
$$
Summing all up, we get 
$$
(\partial_{s}g+\mathcal{L}_{k}^{\varepsilon} g)+(\partial_{s}f+\mathcal{L}_{k}^{\varepsilon} f)\leq0;
$$
this contradiction means that the local maximum can be attained only
at $S=t$. Now we proceed to the equality. Setting for some $\varepsilon,\delta>0$
$$
g_{\varepsilon,\delta}(y_{1},s):=g(y_{1},s)-\delta(t-s)-\varepsilon e^{-s}|y_{1}|^{2},\, f_{\varepsilon,\delta}(y_{2},s):=f(y_{2},s)-\delta(t-s)-\varepsilon e^{-s}|y_{2}|^{2}
$$
and computing $\partial_{s}+\mathcal{L}_{k}^{\varepsilon}$, we come to the previous case for
$\varepsilon,\delta$ small enough (since all the coefficients of the differential
operator are bounded). Passing to the limit as $\varepsilon,\delta\rightarrow0$,
we come to the required assertion.

\textbf{Step 5. Deriving the estimate-1.} Plugging solutions of \eqref{eq:backw}
into identity \eqref{r2}, we get 
$$
\int\phi d\mu_{t}-\int g(x,0)d\mu_{0}=-\int_{0}^{t}\int(A_{k}^{\varepsilon}(x,s)-B_{\mu})\cdot\nabla g(x,s)d\mu_{s}ds,
$$
$$
\int\psi d\sigma_{t}-\int f(x,0)d\sigma_{0}=-\int_{0}^{t}\int(A_{k}^{\varepsilon}(x,s)-B_{\sigma})\cdot\nabla f(x,s)d\sigma_{s}ds.
$$
 Because of \eqref{eq:max_pr_pr_1}, 
 $$
\int\phi d\mu_{t}-\int g(x,0)d\mu_{0}\leq l\int_{0}^{t}\int|A_{k}^{\varepsilon}-B_{\mu}|d\mu_{s}ds.
$$
 Note that 
 \begin{multline*}
\int\psi d\sigma_{t}-\int f(x,0)d\sigma_{0}\leq\int_{0}^{t}\int|A_{k}^{\varepsilon}-B_{\mu}|\cdot|\nabla f|d\sigma_{s}ds
+\int_{0}^{t}\int|B_{\mu}-B_{\sigma}|\cdot|\nabla f|d\sigma_{s}ds\overset{\eqref{eq:max_pr_pr}}{\leq}\\
\leq l\int_{0}^{t}\int|A_{k}^{\varepsilon}-B_{\mu}|d\sigma_{s}ds+
\int_{0}^{t}\int|B_{\mu}-B_{\sigma}|\cdot|\nabla f|d\sigma_{s}ds.
\end{multline*}
 Summing up these inequalities, we get 
 $$
\int\phi d\mu_{t}+\int\psi d\sigma_{t}\leq\int g(x,0)d\mu_{0}+
\int f(x,0)d\sigma_{0}+l\cdot R_{k}^{\varepsilon}+
\int_{0}^{t}\int|B_{\mu}-B_{\sigma}|\cdot|\nabla f|d\sigma_{s}ds,
$$
where 
$$
R_{k}^{\varepsilon}:=\int_{0}^{t}\int|A_{k}^{\varepsilon}-B_{\mu}|d(\mu_{s}+\sigma_{s})ds.
$$
 By virtue of step 4 we have $g(x,0)+f(y,0)\leq h(|x-y|)$.
Thus 
\begin{equation}
\int g(x,0)d\mu_{0}+\int f(x,0)d\sigma_{0}\leq C_h (\mu_0,\sigma_0).
\label{eq:pr-m-m}\end{equation}
 So we get 
 \begin{equation}
\int\phi d\mu_{t}+\int\psi d\sigma_{t}\leq
C_h (\mu_0,\sigma_0)+
l\cdot R_{k}^{\varepsilon}+\int_{0}^{t}\int|B_{\mu}-B_{\sigma}|\cdot|\nabla f|d\sigma_{s}ds.
\label{eq:basic-estimate}\end{equation}
\textbf{ Step 6. Integral bound for $\nabla f$.} The last term in
the right-hand side of \eqref{eq:basic-estimate} is dominated by
$$
\sqrt{\int_{0}^{t}\int\nu^{-1}|B_{\mu}-B_{\sigma}|^{2}d\sigma_{s}ds}\cdot
\sqrt{\int_{0}^{t}\int|\sqrt{Q}\nabla f|^{2}d\sigma_{s}ds}.
$$
To estimate the second multiplier, i.e. 
${\displaystyle \Vert \sqrt{Q}\nabla f\Vert }_{L^{2}(\mathbb{R}^{d}\times[0,T];d\sigma_{s}ds)}^{2}$, 
note that $f^{2}$ is a function of the class $C_{b}^{2,1}(\mathbb{R}^{d}\times[0,t))\cap C(\mathbb{R}^{d}\times[0,T])$
and it can be plugged into the identity \eqref{r2} for the measure
$\sigma$:
\begin{multline*}
\int\psi^{2}d\sigma_{t}-\int f^{2}(x,0)d\sigma_{0}=\int_{0}^{t}\int(\partial_{s}+L_{\sigma})f^{2}d\sigma_{s}ds\\
=\int_{0}^{t}\int2f(\partial_{s}f+\mbox{trace}(QD^{2}f)+\langle B_{\sigma},\nabla f\rangle )+2|\nabla f|^{2}d\sigma_{s}ds\\
=-\int_{0}^{t}\int2f\langle A_{k}^{\varepsilon}-B_{\sigma},\nabla f\rangle +2|\nabla f|^{2}d\sigma_{s}ds.
\end{multline*}
 Hence 
$$
2\int_{0}^{t}\int|\nabla f|^{2}d\sigma_{s}ds\leq\int\psi^{2}d\sigma_{t}-\int f^{2}(x,0)d\sigma_{0}
+2\max|f(x,t)|\int_{0}^{t}\int\nu^{-1/2}|A_{k}^{\varepsilon}-B_{\sigma}|\cdot|\sqrt{Q}\nabla f|d\sigma_{s}ds.
$$
 The maximum principle \eqref{eq:max_pr} and definition \eqref{eq:phib}
imply $\max|f(x,t)|\leq\max|\psi(x)|\leq\Vert h\Vert _{\infty}$.
Taking into account the Cauchy inequality $ab\leq2^{-1}\gamma a^{2}+(2\gamma)^{-1}b^{2}$
with $\gamma=\Vert h\Vert _{\infty}$, we come to 
$$
2\int_{0}^{t}\int|\nabla f|^{2}d\sigma_{s}ds\leq\Vert h\Vert _{\infty}^{2}+
\Vert h\Vert _{\infty}^{2}\nu^{-1}\int_{0}^{t}\int|A_{k}^{\varepsilon}-B_{\sigma}|^{2}d\sigma_{s}ds+
\int_{0}^{t}\int|\sqrt{Q}\nabla f|^{2}d\sigma_{s}ds.
$$
 Cancelling alike terms and recalling \eqref{eq:basic-estimate},
we get 
\begin{equation}
\int\phi d\mu_{t}+\int\psi d\sigma_{t}\leq C_h (\mu_0,\sigma_0)+l\cdot R_{k}^{\varepsilon}+\Vert h\Vert _{\infty}\nu^{-1/2}\cdot r_{k}^{\varepsilon}\cdot\sqrt{\int_{0}^{t}\int|B_{\mu}-B_{\sigma}|^{2}d\sigma_{s}ds},
\label{eq:basic_K}\end{equation}
 where 
 $$
r_{k}^{\varepsilon}:=\sqrt{1+\nu^{-1}\int_{0}^{t}\int|A_{k}^{\varepsilon}-B_{\sigma}|^{2}d\sigma_{s}ds}.
$$
\textbf{ Step 7. Limits as $\varepsilon\rightarrow0$ and $k\rightarrow\infty$.
Deriving the estimate-2.} First of all, we recall that 
$$
A_{k}^{\varepsilon}(x,t)\rightarrow A_{k}(x,t)\,\,\mbox{for a.e. }(x,t)\in\mathbb{R}^{d}\times[0,T],
$$
and the measure $d\sigma_{s}ds$ and $d(\mu_{s}+\sigma_{s})ds$ have
strictly positive densities on $\mathbb{R}^{d}\times[0,T]$ with respect to  the Lebesgue
measure. Thus 
$$
A_{k}^{\varepsilon}(x,t)\rightarrow A_{k}(x,t)\quad d\sigma_{s}ds\,\mbox{-a.e. and }d(\mu_{s}+
\sigma_{s})ds\,\mbox{-a.e.}
$$
Since for each $k$ the mappings $A_{k}^{\varepsilon}$ and $A_{k}$ are
bounded, Lebesgue's dominated convergence theorem yields 
\begin{multline*}
\int_{0}^{t}\int|A_{k}^{\varepsilon}-B_{\sigma}|^{2}d\sigma_{s}ds\rightarrow\int_{0}^{t}\int|A_{k}-B_{\sigma}|^{2}d\sigma_{s}ds,\quad\varepsilon\rightarrow0,\\
\int_{0}^{t}\int|A_{k}^{\varepsilon}-B_{\mu}|d(\mu_{s}+\sigma_{s})ds\rightarrow\int_{0}^{t}\int|A_{k}-B_{\mu}|d(\mu_{s}+\sigma_{s})ds,\quad\varepsilon\rightarrow0.
\end{multline*}
Next, recall that 
$$
\lim_{k\rightarrow\infty}A_{k}(x,t)=B_{\mu}(x,t)\,\,\mbox{for a.e.}\,\,(x,t)\in\mathbb{R}^{d}\times[0,T].
$$
Similarly, taking into account \eqref{eq:B-reg-1},
one can apply the Lebesgue's dominated convergence theorem and get
$$
\lim_{k\rightarrow\infty}\lim_{\varepsilon\rightarrow0}R_{k}^{\varepsilon}=0,
\quad\lim_{k\rightarrow\infty}\lim_{\varepsilon\rightarrow0}r_{k}^{\varepsilon}=
\sqrt{1+\int_{0}^{t}\int\nu^{-1}\cdot|B_{\mu}-B_{\sigma}|^{2}d\sigma_{s}ds}.
$$
 Hence one can pass in \eqref{eq:basic_K} to limits as $\varepsilon\rightarrow0$,
then $k\rightarrow\infty$ and get 
\begin{multline}
\int\phi d\mu_{t}+\int\psi d\sigma_{t}\leq C_h (\mu_0,\sigma_0)\\
+\Vert h\Vert _{\infty}\sqrt{\int_{0}^{t}\int\nu^{-1}|B_{\mu}-B_{\sigma}|^{2}d\sigma_{s}ds}
\cdot\sqrt{1+\nu^{-1}\int_{0}^{t}\int|B_{\mu}-B_{\sigma}|^{2}d\sigma_{s}ds}.
\label{eq:basic_ed}\end{multline}
 Passing to supremum over $(\phi,\psi)\in\Phi_{h}^{b}$
and using step 3, we obtain 
$$
C_h (\mu_t,\sigma_t) \leq C_h (\mu_0,\sigma_0) 
+\Vert h\Vert _{\infty}\sqrt{\int_{0}^{t}\int\nu^{-1}|B_{\mu}-B_{\sigma}|^{2}d\sigma_{s}ds}\cdot
\sqrt{1+\nu^{-1}\int_{0}^{t}\int|B_{\mu}-B_{\sigma}|^{2}d\sigma_{s}ds},
$$
 i.e. the estimate \eqref{eq:basic_est} with $\lambda=0$. \hfill{}$\blacksquare$

\section{Applications to nonlinear equations }

In this section we focus on the applications of the obtained estimate
to the study of the well-posedness of the Cauchy problem for the nonlinear
FPK equations. 

Given a continuous positive function $\alpha$ on $[0,T]$, $\tau\in(0,T]$
and a non-negative continuous function $V(x)$ on $\mathbb{R}^{d}$ with $V(x)\rightarrow+\infty$
as $|x|\rightarrow+\infty$, define classes of measures 
\begin{multline*}
M_{\tau,\alpha}(V)=\{ \mu=(\mu_{t})_{t\in[0,\tau]}:\int V(x)d\mu_{t}\leq\alpha(t),\, t\in[0,\tau]\} ,\\
\quad M_{\tau}(V)=\{ \mu=(\mu_{t})_{t\in[0,\tau]}:\sup_{t\in[0,\tau]}\int V(x)d\mu_{t}<+\infty\} .
\end{multline*}
Throughout the section we assume that a non-degenerate $d\times d$-matrix
$Q(x,t)$ satisfying ${\rm (A1)}$ is fixed. Suppose that for each
measure $\mu=\mu_{t}dt\in M_{T}(V)$ a Borel mapping 
$$
B(\mu,\cdot,\cdot)\equiv B(\mu):\mathbb{R}^{d}\times[0,T]\rightarrow\mathbb{R}^{d}
$$
is defined. Consider the Cauchy problem for a nonlinear FPK equation
\begin{equation}
\partial_{t}\mu_{t}=\mbox{trace}(Q(x,t)D^{2}\mu_{t})-\mbox{div}(B(\mu,x,t)\mu_{t}),\quad\mu_{t}|_{t=0}=\mu_{0}.
\label{eq:nonl_label}\end{equation}
Again denote the elements of the diffusion matrix $Q(x,t)$ by $q^{ij}(x,t),\,1\leq i,j\leq d$
and the elements of the vector drift $B(\mu,x,t)$ by $b^{j}(\mu,x,t),\,1\leq j\leq d$.
Set 
$$
L_{\mu}\phi=q^{ij}(x,t)\partial_{x_{i}x_{j}}^{2}\phi+b^{i}(\mu,x,t)\partial_{x_{i}}\phi,
$$
where summation over all repeated indices is taken. As earlier, we
call the measure $\mu=\mu_{t}dt,\, t\in[0,T]$ a solution
to \eqref{eq:nonl_label}, if the identity \eqref{r1} holds with
$L_{\mu}$instead of $L$. Introduce the following assumptions on
the drift: 

(B1) The drift term $B$ is $\lambda$-dissipative in $x$, i.e. for
every measure $\mu\in M_{T}(V)$ 
\begin{equation}
\langle B(\mu,x,t)-B(\mu,y,t),x-y\rangle _{\mathbb{R}^{d}}\leq\lambda\Vert x-y\Vert ^{2}
\label{eq:B-1-1}\end{equation}
for all $x,y\in\mathbb{R}^{d}$ and all $t\in[0,T]$.

(B2)  for all measures $\mu$ and $\sigma$ from $M_{T}(V)$
\begin{equation}
B(\mu,x,t)-\lambda x\,\in L^{2}(\mathbb{R}^{d}\times[0,T],d(\mu_{s}+\sigma_{s})ds).
\label{eq:B-reg-1-1}\end{equation}

We start with the question of uniqueness and stability of the probability
solution to \eqref{eq:nonl_label}. As earlier, we assume that some
continuous non-decreasing monotone bounded cost function  $h$ with
$h(0)=0$ is fixed. Given a non-negative non-decreasing
function $G$, denote 
$$
G^{*}(r):={\displaystyle \int_{r}^{1}\frac{du}{G^{2}(\sqrt{u})}}.
$$
 \begin{corollary} \label{uni+est} Fix some non-negative continuous
function $V(x)$ on $\mathbb{R}^{d}$ with $V(x)\rightarrow+\infty$ as $|x|\rightarrow+\infty$
such that $V\in L^{1}(\mathbb{R}^{d};\mu_{0})\cap L^{1}(\mathbb{R}^{d};\sigma_{0})$. Assume
that the coefficients of the equation \eqref{eq:nonl_label} satisfy
${\rm (A1)},$ ${\rm (B1)}$ and ${\rm (B2)}$ with this $V$. Moreover,
assume that that for each two measures $\mu=(\mu_{t})_{t\in[0,T]}$
and $(\sigma_{t})_{t\in[0,T]}$ from $M_{T}(V)$ 
\begin{equation}
|B(\mu,x,t)-B(\sigma,x,t)|\leq\sqrt{V(x)}G(C_h (\mu_t,\sigma_t)
\label{eq:lip}\end{equation}
for some non-negative increasing function $G$ such that $G^{*}(0)=+\infty$. 

Then each two solutions $(\mu_{t})_{t\in[0,T]}$ and $(\sigma_{t})_{t\in[0,T]}$
of the problem \eqref{eq:nonl_label} from the class $M_{T}(V)$ with
initial data $\mu_{0}$ and $\sigma_{0}$ respectively satisfy 
$$
C_{h_{\lambda t}} (\mu_t,\sigma_t)\leq\Bigl((G^{*})^{-1}\Bigl
(G^{*}(2(C_h (\mu_0,\sigma_0))^{2})-ct\Bigr)\Bigr){}^{1/2}
$$
for all $t\in[0,T]$; here $(G^{*})^{-1}$ is the inverse to $G^{*}$
function, and $c>0$ is some positive constant. \end{corollary}

\begin{example} Assumptions ${\rm (B1)}$, ${\rm (B2)}$ and \eqref{eq:lip}
are fulfilled, for example, for drift terms of the form 
$$
B(\mu,x,t)=H(x)\int k(x,y)d\mu_{t}(y)
$$
with $0\leq H(x)\leq\sqrt{V(x)}$ and a $\lambda$- dissipative in
the first variable kernel $k(\cdot,\cdot)$ such that 
$$
|k(x,y)-k(z,y)|\leq h(|x-y|).
$$

\end{example}

\textbf{Proof of Corollary \ref{uni+est}.} First of all, if $\mu$
is a solution to \eqref{eq:nonl_label} and assumptions of Corollary
\ref{uni+est} are fulfilled, then the linear FPK equation 
$$
\partial_{t}\rho_{t}=\mbox{trace}(Q(x,t)D^{2}\rho_{t})-
\mbox{div}(B(\mu,x,t)\rho_{t}),\quad\rho_{t}|_{t=0}=\mu_{0}
$$
has a solution $\rho=\mu$ and it satisfies assumptions of Theorem
\ref{est}; similarly does $\sigma$. Hence one can apply \eqref{eq:basic_est}
with $B_{\mu}(\cdot,\cdot)=B(\mu,\cdot,\cdot)$
and $B_{\sigma}(\cdot,\cdot)=B(\sigma,\cdot,\cdot)$. 

Next, arguing as on Step 1 of the proof of Theorem \ref{est}, we
can assume that the drift term $B$ is dissipative. With condition
\eqref{eq:lip} in hand, the estimate \eqref{eq:basic_est} takes
the form 
\begin{equation}
C_h (\mu_t,\sigma_t) \leq C_h (\mu_0,\sigma_0)
+\Vert h\Vert _{\infty}\sqrt{\nu^{-1}a}
\sqrt{\int_{0}^{t}G^{2}(C_h (\mu_s,\sigma_s))ds}
\cdot\sqrt{1+\nu^{-1}a\int_{0}^{t}G^{2}(C_h (\mu_s,\sigma_s))ds},
\label{eq:inequality}\end{equation}
where $a=\sup_{t\in[0,T]}\int V(x)d\mu_{t}<+\infty$ and $\nu$ is
the ellipticity constant of $Q$. Note that 
$C_h (\mu_t,\sigma_t)\leq\Vert h\Vert _{\infty}$. 
Then
\eqref{eq:inequality} can be reduced to a weaker inequality 
\begin{equation}
C_h (\mu_t,\sigma_t)\leq C_h (\mu_0,\sigma_0)+K\sqrt{\int_{0}^{t}G^{2}(C_h (\mu_s,\sigma_s))ds}
\label{eq:int_ineq}\end{equation}
 with $K=\Vert h\Vert _{\infty}\sqrt{\nu^{-1}\alpha}\cdot\sqrt{1+\nu^{-1}\cdot TG^{2}(\Vert h\Vert _{\infty})}.$
Squaring \eqref{eq:int_ineq} and using the inequality $(b+c)^{2}\leq2b^{2}+2c^{2}$,
we get 
$$
C_h (\mu_t,\sigma_t)^{2}\leq2C_h (\mu_0,\sigma_0)^{2}+2K^{2}\int_{0}^{t}G^{2}(C_h (\mu_s,\sigma_s))ds.
$$
 If $\mu_{0}=\sigma_{0}$, then uniqueness follows immediately due
to explicit integration. In the general case the Gronwall type inequality
(for example, \cite[Th. 27]{gronwall-type}) implies 
$$
C_h (\mu_t,\sigma_t) \leq((G^{*})^{-1}\Bigl(G^{*}(2(C_h (\mu_0,\sigma_0))^{2})-2K^{2}t\Bigr))^{1/2}.
$$
\hfill{}$\square$

A particular special case $G(u)=u$ of this latter estimate is especially
interesting:

\begin{corollary}\label{znak_l} Let $\mu$ and $\sigma$ be two
solutions to \eqref{eq:nonl_label} as in Theorem \ref{uni+est} with
$G(u)=u$. Then for some $N>0$ 
$$
C_{h_{\lambda t}} (\mu_t,\sigma_t) \leq \sqrt{2}C_{h} (\mu_0,\sigma_0) e^{Nt}.
$$
 In particular, if the drift is dissipative ($\lambda=0$) or $\lambda<0$,
then 
$$
C_h (\mu_t,\sigma_t) \leq\sqrt{2}C_{h} (\mu_0,\sigma_0) e^{Nt}.
$$
 \end{corollary}

In some cases the estimate \eqref{eq:basic_est}  enables to establish
existence of a solution to the nonlinear equation \eqref{eq:nonl_label}.
To show this, consider $h(r)=\min\{ |r|^{p},1\} $
for some $p\geq1$. Recall that in this case $C_{h}^{1/p}(\mu_{t},\sigma_{t})$ is
a metric and turns the space of probability measures into a complete
metric space. Moreover, convergence with respect to this metric is
equivalent to weak convergence (see \cite[Th. 1.1.9]{obzor_Monge}). 

\begin{corollary} \label{ex_schauder} Suppose there exists a function
$V\in C^{2}(\mathbb{R}^{d})$, $V\geq1$ such that $V(x)\rightarrow+\infty$
as $|x|\rightarrow+\infty$ and there exists positive function $\Lambda$
on $[0,+\infty)$ such that 
$$
(L_{\mu}V)(x,t)\leq\Lambda(\alpha(t))(1+V(x))
$$
for each $\alpha\in C^{+}[0,T]$, $\tau\in[0,T]$, each $(x,t)\in\mathbb{R}^{d}\times[0,T]$
and each $\mu\in M_{\tau,\alpha}(V)$. Assume that the
coefficients in \eqref{eq:nonl_label} satisfy ${\rm (A1)},$ ${\rm (B1)}$
and ${\rm (B2)}$ with this function $V$. Assume that $B(\sigma^{n})\rightarrow B(\sigma)$
in $L^{2}(\mathbb{R}^{d}\times[0,T],d\sigma_{s}ds)$ as
$n\rightarrow\infty$ if measures $\sigma^{n}(dx\, dt)=\sigma_{t}^{n}(dx)dt$
weakly converge to a measure $\sigma(dx\, dt)=\sigma_{t}(dx)dt$ on
the strip $\mathbb{R}^{d}\times[0,T]$. 

Then for every probability measure $\mu^{*}$ such that $V\in L^{1}(\mathbb{R}^{d};\mu^{*})$,
there exists a (local) probability solution $\mu=(\mu_{t})_{t\in[0,\tau]}$
to \eqref{eq:nonl_label} with initial condition $\mu^{*}$. \end{corollary} 

\begin{example} Let $k(x,y)$ be a bounded function, $\lambda$-
dissipative in the first variable for every $y\in\mathbb{R}^{d}$. Let $Q(x,t)$
be a matrix satisfying ${\rm (A1)}$. Then the Cauchy problem \eqref{eq:nonl_label}
with 
$$
B(\mu,x,t)=\int k(x,y)d\mu_{t}(y)
$$
satisfies all assumptions of Theorem \ref{ex_schauder} with $V(x)=1+|x|^{2}$
and any probability measure $\nu$ with finite second moment. \end{example}

\begin{example} Let $V>0$ be some $C^{2}$-function on $\mathbb{R}^{d}$ with
at least linear growth. Let $g(x)$ be a $\lambda$- dissipative function
such that $|g|\leq\sqrt{V}$. Let $Q(x,t)$ be a matrix satisfying
${\rm (A1)}$. Then the Cauchy problem \eqref{eq:nonl_label} with
$$
B(\mu,x,t)=g(x)\int k(y)d\mu_{t}(y)
$$
with some non-negative continuous bounded kernel $k(y)$ satisfies
all assumptions of Theorem \ref{ex_schauder} with any probability
measure $\nu$ that integrates $V$. \end{example}

\begin{example} Fix $\alpha\in(0,1)$ and a matrix $Q$ satisfying
${\rm (A1)}$. Then the Cauchy problem \eqref{eq:nonl_label} with
$$
B(\mu,x,t)=-(|x|^{\alpha-1}x)*\mu_{t}
$$
satisfies all assumptions of Theorem \ref{ex_schauder} with $V(x)=1+|x|^{2}$
and any probability measure $\nu$ with finite second moment (cf.
\cite[Proposition 2.1]{ManShap}). \end{example}

\textbf{Proof.} As earlier, without loss of generality, the drift
term is dissipative. Let $\sigma\in M_{\tau,\alpha}(V)$ for some
$\tau,\alpha$. Consider 
$$
\partial_{t}\mu_{t}=\mbox{trace}(Q(x,t)D^{2}\mu_{t})-\mbox{div}(B(\sigma,x,t)\mu_{t}),\quad\mu_{0}=\mu^{*}.
$$
Note that the dissipativity of the drift ensures that it is bounded  localy
in $(x,t)$. Hence under the assumptions of the theorem there
exists a unique probability solution $\mu=(\mu_{t})_{t\in[0,\tau]}$
in $M_{\tau}(V)$ (see \cite[Theorem 3.6]{ManShapJDDE}). Therefore
the mapping $\Theta:\, M_{\tau,\alpha}(V)\rightarrow M_{\tau}(V)$
$$
\mu=\Theta(\sigma)\Longleftrightarrow\partial_{t}\mu_{t}=\mbox{trace}(Q(x,t)D^{2}\mu_{t})-
\mbox{div}(B(\sigma,x,t)\mu_{t}),\quad\mu_{0}=\mu^{*}
$$
is correctly defined. It is obvious that the solutions to \eqref{eq:nonl_label}
are exactly the fixed points of the mapping $\Theta$.

Define subclass $N_{\tau,\alpha}(V)$ of the class $M_{\tau,\alpha}(V)$
as follows: 
$$
N_{\tau,\alpha}(V):=\{ \mu\in M_{\tau,\alpha}(V):
|\int\varphi(x)d(\mu_{t}-\mu_{s})|\leq K(\tau,\alpha,\varphi)\cdot|t-s|\,\,
\forall\varphi\in C_{0}^{\infty}(\mathbb{R}^{d})\} ,
$$
where 
$$
K(\tau,\alpha,\varphi):=\sup\{ |L_{\mu}\varphi(x,t)|,\,(x,t)\in\mathbb{R}^{d}\times[0,\tau],\,\mu\in M_{\tau,\alpha}(V)\}.
$$
Obviously $N_{\tau,\alpha}$ is a convex set. By virtue of \cite[Corollary 4]{ManShap}
there exist $\bar{\alpha}(t)>0$ and $\bar{\tau}\in(0,T]$ such that
$\Theta(N_{\bar{\tau},\bar{\alpha}}(V))\subset N_{\bar{\tau},\bar{\alpha}}(V)$.
Moreover, the class $N_{\bar{\tau},\bar{\alpha}}(V)$ is a compact
set in the topology of weak convergence of measures on the strip $\mathbb{R}^{d}\times[0,\tau]$
by \cite[Corollary 1]{ManShap}. Let us check that continuity of the  mapping $\Theta$
 on $N_{\bar{\tau},\bar{\alpha}}(V)$. Suppose that the
sequence $\sigma^{n}=(\sigma_{t}^{n})\in N_{\bar{\tau},\bar{\alpha}}(V)$
weakly converges to $\sigma=(\sigma_{t})\in N_{\bar{\tau},\bar{\alpha}}(V)$.
Set $\mu^{n}:=\Theta(\sigma^{n})$, $\mu:=\Theta(\sigma)$.
Due to \eqref{eq:basic-estimate} we have 
$$
C_{h}(\mu_{t}^{n},\mu_{t})\leq\sqrt{\int_{0}^{\bar{\tau}}\int|B(\sigma^{n})-B(\sigma)|^{2}d\sigma_{s}ds}\cdot\sqrt{1+\int_{0}^{\bar{\tau}}\int|B(\sigma^{n})-B(\sigma)|^{2}d\sigma_{s}ds}.
$$
Our conditions imply that the right-hand side goes to zero as $n\rightarrow\infty$.
Hence $\mu_{t}^{n}$ converges to $\mu_{t}$ with respect to the metric
$C_{h}^{1/p}$ and thus converges weakly. Let us show that $\mu^{n}$
converges to $\mu$ on the strip $\mathbb{R}^{d}\times[0,\bar{\tau}]$. Fix some
continuous bounded function $\zeta(x,t)$. Then for each $t\in[0,\bar{\tau}]$
we have 
$$
\int\zeta(x,t)d\mu_{t}^{n}\rightarrow\int\zeta(x,t)d\mu_{t},\,\, n\rightarrow\infty.
$$
Since the measures $\mu_{t}^{n}$ are probability measures and $\zeta$
is bounded, the integrals on the right-hand side are uniformly bounded
and pointwise (with respect to $t\in[0,\bar{\tau}]$) convergent to ${\displaystyle \int\zeta(x,t)d\mu_{t}}$.
Therefore Lebesgue's dominated convergence theorem ensures 
$$
\int_{0}^{\bar{\tau}}\int\zeta(x,t)d\mu_{t}^{n}dt\rightarrow\int_{0}^{\bar{\tau}}\int\zeta(x,t)d\mu_{t}dt,
\,\, n\rightarrow\infty.
$$
By definition this means that the sequence $\mu^{n}$ converges weakly
to $\mu$ on the strip $\mathbb{R}^{d}\times[0,\bar{\tau}]$. 

Summarizing, we have a continuous mapping $\Theta$ on a convex compact
set $N_{\bar{\tau},\bar{\alpha}}(V)$ and maps it onto itself. The
Schauder fixed-point theorem ensures that there exists a fixed point of
$\Theta$ in $N_{\bar{\tau},\bar{\alpha}}(V)$, i.e. there exists
a solution $\mu=(\mu)_{t\in[0,\bar{\tau}]}$
to \eqref{eq:nonl_label} with initial condition $\mu^{*}$. \hfill{}$\square$

~

\textbf{Acknowledgements.}

The work was partially supported by RFBR grants 14-01-00237 and 15-31-20082.
The author is grateful to M. Zaal for fruitful discussions and
to A.I. Nazarov for valuable comments and remarks.

\end{document}